\documentclass{article}

\usepackage{dsfont}
\usepackage{fixmath}
\usepackage{hyperref}
\usepackage{amsmath}
\usepackage{amsthm}
\usepackage{amssymb}
\usepackage{mathtools}
\usepackage{mathabx,epsfig}
\usepackage{tikz-cd}

\newtheorem{theorem}{Theorem} 

\newtheorem{definition}[theorem]{Definition}

\newtheorem{lemma}[theorem]{Lemma}

\newtheorem{remark}[theorem]{Remark}

\newtheorem{corollary}[theorem]{Corollary}

\newtheorem{example}[theorem]{Example}

\newtheorem{question}[theorem]{Question}

\newtheorem{conjecture}[theorem]{Conjecture}

\def\acts{\mathrel{\reflectbox{$\righttoleftarrow$}}}

\newcommand{\Z}{\mathbb{Z}}
\newcommand{\C}{\mathbb{C}}
\newcommand{\N}{\mathbb{N}}

\title{Dynamical Perturbing and $C^*$-algebra Lifting Problems}
\author{Samantha Pilgrim}
\date{\today}

\begin{document}
\maketitle

\section{Introduction}

It is often interesting to consider whether set maps which are ``close to being morphisms" (in some sense) are ``close to morphisms" (in some other sense).  In the category of group representations, for instance, one can ask, for a sequence of set maps maps $\rho_n:\Gamma\to \mathcal{B}(V_n)$ with $\|\rho_n(gh)-\rho_n(g)\rho_n(h)\|\to 0$ for all $g, h\in \Gamma$, whether there exist representations $\tilde{\rho}_n:\Gamma\to \mathcal{B}(V_n)$ with $\|\rho_n(g) - \tilde{\rho}_n(g)\|\to 0$.  In the category of $C^*$-algebras, this property is related to weak semiprojectivity of the group $C^*$-algebra $C^*\Gamma$.  Almost commuting unitaries (asymptotic representations of $\mathbb{Z}^2$) in particular have seen significant study, and are known to encode homological information about $C^*$-algebras \cite{Bott_almost_commuting}.  In the context of topology (commutative $C^*$-algebras), semiprojectivity coincides with the notion of ANR \cite{semiprojectivity_commutative}.  

	In a very different context, one has the notion of flexible stability in permutations \cite{Lazarovich2019SurfaceGA}, wherein asymptotic homomorphisms into finite symmetric groups (in the sense of the normalized hamming distance) can be perturbed to actual homomorphisms.  Any sofic group which is flexibly stable in permutations is then residually finite.  Similarly, an MF $C^*$-algebra which is semiprojective is residually finite dimensional \cite[Corollary 2.14]{BlackadarSimple}.  These stability properties then have two significant flavors of applications: on the one hand they can be used to infer stronger finite (dimensional) approximation properties from weaker ones, and on the other they can potentially be used to verify that even weak approximations do not exist for certain examples (e.g. potentially proving an example of a non-sofic group by proving a group is flexibly stable in permutations but not residually finite \cite{MR3900769}).  
	
This note initiates study of an analogous question in the category of topological group actions.  More precisely, we ask the following:
	
\question{Suppose $\Gamma$ is a group, $(X_n, d_n)$ are compact metric spaces, and $\alpha_n: \Gamma\to \text{Homeo}(X_n)$ are set maps such that

$$\sup_{x\in X_n} d_n(\alpha_n(\gamma)\circ \alpha_n(\delta)(x), \alpha_n(\gamma\delta)(x))\to 0$$

\noindent as $n\to\infty$ for all $\gamma, \delta\in \Gamma$. Are there homomorphisms $\widetilde{\alpha}_n:\Gamma\to \text{Homeo}(X_n)$ such that $\sup_{x\in X_n} d_n(\alpha_n(\gamma)(x), \widetilde{\alpha}_n(\gamma)(x))\to 0$ for all $\gamma\in \Gamma$?}
\normalfont
	
We will also consider two special cases of this question in particular. 
	
\begin{description}
	\item(1) Given a sequence of set maps $\alpha_n: \Gamma\to \text{Homeo}(X)$ such that $\sup_{x\in X}d_X(\alpha_n(\gamma\delta)(x), \alpha_n(\gamma)\circ \alpha_n(\delta)(x))\to 0$ for all $\gamma, \delta\in \Gamma$, is there a homomorphism $\tilde{\alpha}: \Gamma\to \text{Homeo}(X)$ such that $\sup_{x\in X} d_X(\tilde{\alpha}(\gamma)(x), \alpha_n(\gamma)(x))\to 0$ for all $\gamma\in \Gamma$?
	\item(2) Given finite sets $E_n\subset X$ and a sequence of set maps $\alpha_n: \Gamma\to \text{Sym}(E_n)$ such that $\max_{e\in E_n}d_X(\alpha_n(\gamma\delta)(e), \alpha_n(\gamma)\circ \alpha_n(\delta)(e))\to 0$ for all $\gamma, \delta\in \Gamma$, are there homomorphisms $\tilde{\alpha}_n: \Gamma\to \text{Sym}(E_n)$ such that $\max_{e\in E_n}d_X(\tilde{\alpha}_n(\gamma)(e), \alpha_n(\gamma)(e))\to 0$ for all $\gamma\in \Gamma$?	
\end{description}

\normalfont \noindent We call the sequence $(\alpha_n)$ an almost-action and the problem of finding an honest action nearby the dynamical perturbing problem associated to $(\alpha_n)$.  Problems of type (2) above are of interest to us because current results \cite[Lemma 5.1]{kerr2011residually} can already produce many examples of such almost-actions. One can also distinguish further subclasses of these problems by prescribing whether or not the convergence is uniform over elements of $\Gamma$ (see \ref{main definition}).  

	We show that, under certain assumptions, these perturbing problems always have solutions.  For problems of type (1), we will see in the last section that they always have solutions when $X$ is a Cantor set and $\Gamma$ is finite \ref{pretty Cantor sets theorem}, motivating us to conjecture that this is the case for any compact space $X$.  For (2), we show \ref{cantor sets theorem} the answer is always yes when $\Gamma$ is virtually free, $X$ is a Cantor set, and the almost actions come from \cite[Lemma 5.1]{kerr2011residually}.  This allows us to show actions of virtually free groups on Cantor sets which fix a measure of full support are residually finite \ref{main corollary} in the sense of \cite{kerr2011residually}, which implies, essentially, that their reduced crossed product $C^*$-algebras have finite-dimensional approximation properties as nice as those of the acting group's reduced $C^*$-algebra.  
	
	In an effort to obtain further positive results through functional-analytic methods, we also considered applying the theory of semiprojective $C^*$-algebras, that is, algebras whose approximate representations are close to honest representations.  This property appears at first blush to be just what we need, until one realizes that almost-actions of the kind considered above only give rise to operators in $M_n(\C)$ which are approximately a representation in a pointwise sense as opposed to in norm (as the estimate depends on continuity of a function).  This problem might be surmountable when $\Gamma$ is finite, but even then it's not obvious.  Still, this approach motivated us to investigate a natural, conditional version of semiprojectivity.  We prove that morphisms of finite-dimensional $C^*$-algebras have this property \ref{finite dimensional algebras cond semiprojective} and, in doing so, recapitulate some of the results from \cite{Blackadar1985} with perhaps slightly more digestible proofs.  
	
	Although these topological perturbing problems do not appear to relate to semiprojectivity-type properties of the group or crossed product $C^*$-algebras in the way one might hope, there are still interactions between the two worlds.  Representations constructed from almost-actions can be used to show that even nice inclusions of $1$-dimensional non-commutative CW complexes ($1$-NCCWs) are not conditionally semiprojective \ref{counter-example} in the sense of \ref{conditional semiprojectivity def}.  This is somewhat surprising because, as mentioned, finite-dimensional $C^*$-algebras \textit{are} conditionally semiprojective \cite{+1998+101+143}, and furthermore the semiprojectivity of $1$-NCCWs has been proved precisely by regarding them as special pullbacks of finite-dimensional $C^*$-algebras, and showing that semiprojectivity is preserved when taking such pullbacks \cite{endersfinitecodim}.  
	
	In the last section, we will also see that perturbing problems associated to almost-actions of type (1)  above by finite groups (in fact a slightly larger class) have a purely $C^*$-algebraic characterization in terms of lifting not a homomorphism, but the structure of a Cartan pair \ref{Cartan lift definition}, \ref{Cartan characterization}. A similar result \ref{first Cartan characterization} also holds for almost-actions of type (2) when $\Gamma$ is virtually free and the almost-action comes from \cite[Lemma 5.1]{kerr2011residually}.


	


	Despite the progress documented here, we do not yet have a satisfying account of the boundaries for these types of results.  A possible way forward in this regard would be to consider the almost action $\alpha_n: \Z^2\to \text{Homeo}(T^n)$ given by a specific sequence of almost-commuting unitaries \cite{3d91146b-74b6-3a3d-ba78-139264f362e8}.  In this case, there can be no \textit{isometric} solution to the associated perturbing problem, but it remains to be seen whether there are homeomorphisms of $T^n$ which give a solution.  It is worth noting that the sort of stability we are asking about may not be an exclusively low-dimensional phenomenon as, for example, semiprojectivity of $C^*$-algebras is. It is possible that this notion of stability is significantly more flexible than semiprojectivity (and may, for example, be more closely related to flexible stability, see \cite{Lazarovich2019SurfaceGA}).

\section{Preliminaries}

Throughout, $\Gamma$ will be a countable discrete group, $X$ will be a compact metric space, and $\Gamma\acts X$ will be an action by homeomorphisms.  We denote the action of $\gamma$ on $x\in X$ by $\gamma\cdot x$.  

\subsection{Finite (dimensional) approximations}

\definition{An action $\Gamma\acts X$ is residually finite if for every finite subset $F\subset \Gamma$ and $\epsilon>0$ there is a finite set $E$, a map $\zeta: E\to X$, and an action $\Gamma\acts E$ such that $d_X(\zeta(\gamma\cdot e), \gamma\cdot \zeta(e))<\epsilon$ for all $\gamma\in F$ and $e\in E$.  If $X$ has no isolated points, $\zeta$ can be taken to be an inclusion.}

\definition{An action of $\Gamma$ on a $C^*$-algebra $A$ is quasidiagonal if, for all $\epsilon>0$ and finite sets $S\subset \Gamma$ and $\mathcal{F}\subset A$, there is an $n$, a unital, completely-positive map $\phi: A\to M_n(\C)$, and an action of $\Gamma$ on $M_n(\C)$ such that 

\begin{description}

\item(i) $ \|\phi(a)\| + \epsilon \geq \|a\|$ for all $a\in \mathcal{F}$

\item(ii) $\|\phi(ab) - \phi(a)\phi(b)\|<\epsilon$ for all $a, b\in \mathcal{F}$  

\item(iii) $\|\phi(\gamma\cdot a) - \gamma\cdot \phi(a)\|<\epsilon$ for all $a\in \mathcal{F}$ and $\gamma\in S$.  

\end{description} }

\normalfont It is straightforward to show that a residually finite action $\Gamma\acts X$ is quasidiagonal as an action $\Gamma\acts C(X)$.  One motivation for studying these properties is that they give rise to finite-dimensional approximations of crossed products.  In particular, a crossed product by a quasidiagonal action is quasidiagonal or MF exactly when the $C^*$-algebra of the acting group is quasidiagonal or MF \cite[Theorems 3.4 and 3.5]{kerr2011residually}.  

\subsection{$C^*$-algebra lifting properties}
\normalfont

We will use $\mathcal{A}$, $\mathcal{B}$, and $\mathcal{C}$ for $C^*$-algebras.  If $J = \overline{\cup_n J_n}$ is an increasing union of ideals in $\mathcal{C}$, $\pi: C\to C/J$ will denote the quotient by $J$ and $\pi_n: \mathcal{C}\to \mathcal{C}/J_n$ will denote the quotient by $J_n$.  If $\phi: A\to C/J$ is a homomorphism, we will say a homomorphism $\tilde{\phi}: A\to C/J_n$ lifts $\phi$ if $\pi\circ \tilde{\phi} = \phi$.  

\definition{A $C^*$-algebra $\mathcal{A}$ is called semiprojective if whenever $\rho: A\to \mathcal{C}/J$ is a representation with $J = \cup_n J_n$ an increasing union of ideals, there is a lift $\tilde{\rho}: \mathcal{A}\to \mathcal{C}/J_n$ for some $n$.  We say $\mathcal{A}$ is weakly-semiprojective if this holds only when $\mathcal{C}$ is of the form $\prod_n \mathcal{C}_n$ and $J_n = \bigoplus_{k=1}^n \mathcal{C}_k$, and matrcially-semiprojective if the previous statement holds only when $\mathcal{C}_n = M_{k_n}(\C)$ for all $n$.  


An asymptotic $*$-homomorphism (or asymptotic representation) is a sequence of set maps $\rho_n:\mathcal{A} \to \mathcal{C}_n$ such that

\begin{description}
\item(1) $\|\rho_n(\lambda a+b) - \lambda\rho_n(a) + \rho_n(b)\|\to 0$
\item(2) $\|\rho_n(a^*) - \rho_n(a)^*\|\to 0$
\item(3) $\|\rho_n(ab) - \rho_n(a)\rho_n(b)\|\to 0$
\end{description}
\noindent
as $n\to\infty$ for all $a, b\in \mathcal{A}$ and $\lambda\in \mathbb{C}$.  }

\vspace{3mm}
\normalfont
\noindent The following two lemmas are used in \cite{Blackadar1985} to show that $M_n(\C)$ is semiprojective for all $n$, and that semiprojectivity is preserved when taking finite direct sums (from which it can be concluded that finite dimensional $C^*$-algebras are semiprojective).  We will restate them in a slightly different way and give a few additional details as to their proofs.  

\lemma{Suppose $\rho: \mathcal{A}\to \mathcal{B}/\overline{\cup_n J_n}$ is a $*$-homomorphism.  One can therefore choose an asymptotic $*$-homomorphism $\rho_n: \mathcal{A}\to \mathcal{B}/J_n$ (if $\rho$ is in fact induced by an asymptotic representation, simply choose that one).  Suppose $p\in \mathcal{A}$ is a projection.  Then for any $\epsilon>0$ there is $n$ and a projection $\tilde{p}$ in $B/J_n$ with $\pi(\tilde{p}) = p$ and $\|\rho_n(p) - \tilde{p}\|<\epsilon$.  If $q\in \mathcal{B}/\overline{\cup_n J_n}$ is another such projection orthogonal to $p$, we can arrange that $\tilde{q}$ is orthogonal to $\tilde{q}$.  If $p_1, \ldots, p_k\in \mathcal{A}$ are orthogonal projections such that $p_1+\ldots + p_k = 1_\mathcal{A}$, we can further arrange that $\tilde{p}_1 + \cdots + \tilde{p}_k = 1_{\mathcal{B}/J_n}$.  }

\begin{proof}
This is essentially a rephrasing of \cite[Proposition 2.14]{Blackadar1985}.  The application of functional calculus perturbs $\rho_n(p)$ by $<\epsilon$ if $\rho_n(p)^2 - \rho_n(p)<\epsilon$.  For the last part, we can find lifts $\tilde{p}$ and $\tilde{q}$, which are already close to being orthogonal and then replace $\tilde{q}$ by $(1_{\mathcal{B}/J_n} - \tilde{p})\tilde{q}(1_{\mathcal{B}/J_n} - \tilde{p})$, which doesn't change it by very much.   For the last part, we can find lifts $\tilde{p}_i$ for $1\leq i\leq k-1$, then put $\tilde{p}_k = 1_{\mathcal{B}/J_n} - (\tilde{p}_1 + \cdots + \tilde{p}_{k-1})$, which again doesn't change $p_k$ by very much.  \end{proof}

\lemma{With the same setup as the previous lemma, suppose $v\in \mathcal{A}$ is a partial isometry with $\rho(v^*v) = p$ and $\rho(vv^*) = q$.  Apply the previous lemma to the projections $p$ and $q$ to find projections $\tilde{p}$ and $\tilde{q}$ in $\mathcal{B}/J_n$.  Then for all $\epsilon>0$ there is $n$ and a partial isometry $\tilde{v}\in B/J_n$ such that $\pi(\tilde{v}) = v$, $\tilde{v}\tilde{v}^* = \tilde{p}$, $\tilde{v}^*\tilde{v} = \tilde{q}$, and $\|\rho_n(v) - \tilde{v}\|<\epsilon$.  } \label{lifting partial isometries}

\begin{proof}
Much like the previous lemma, this is a rephrasing of \cite[Proposition 2.23]{Blackadar1985}.  Following through the functional calculus in the proof of this proposition, let $f$ be a continuous function which is $0$ near $0$ and for which $f(\lambda) = \lambda^{-1/2}$ near $1$.  We claim that $\rho_n(v)$ is close to $w = \rho_n(v)f(\rho_n(v)^*\rho_n(v))$.  Since $vv*v=v$ and $(\rho_n)$ is an asymptotic homomorphism, $\|\rho_n(v)-\rho_n(v)\rho_n(v)^*\rho_n(v)\|\to 0$ as $n\to\infty$.  By the same reasoning as the previous lemma, $\sigma(\rho_n(v)^*\rho_n(v))\subset (-\epsilon, \epsilon)\cup (1-\epsilon, 1+\epsilon)$ for any $\epsilon$ for $n$ sufficiently large, and so likewise $\|\rho_n(v)^*\rho_n(v)-f(\rho_n(v)^*\rho_n(v))\|\to 0$.  Hence, $\|\rho_n(v)-\rho_n(v)f(\rho_n(v)^*\rho_n(v))\|\to 0$ as $n\to\infty$.  

%
%

Now, $w$ is a partial isometry with $w^*w = p'_n$ and $ww^* = q'_n$ with $p'_n$ as close as we want to $p_n$ for $n$ large and similarly for $q_n$.  It is shown in \cite[Section 6 Lemma 4]{operatorKfunctor} that if $z_1=(2p'-1)(2\tilde{p} - 1) + 1$ and $z_2 = (2q' - 1)(2\tilde{q} - 1) + 1$, then $z_1$ is invertible with $z_1p'z_1^{-1} = \tilde{p}$ and similarly for $z_2$.  Furthermore, $v_1 = z_1(z_1z_1^*)^{-1/2}$ and $v_2 = z_2(z_2z_2^*)^{-1/2}$ are unitaries and $v_2wv_1$ is the desired partial isometry.  Since each $z_i$ is close to twice the identity, functional calculus shows each $v_i$ is close to the identity.  Hence, $v_2wv_1$ is close to $w$ and therefore close to $\rho_n(v)$ for $n$ sufficiently large.  \end{proof}

\normalfont \noindent We conclude this subsection with some examples of semiprojective crossed products.  Also related to this topic is the notion of equivariant semiprojectivity \cite{PHILLIPS2015929}, but we leave further exposition of the relationship between the two concepts for another time.  

\example{Let $\Gamma$ be finite, $X$ a Cantor set, and $\alpha: \Gamma\to \text{Aut}(X)$.  The crossed product $C(X)\rtimes \Gamma$ is weakly semiprojective.  } \label{weakly semiprojective}

\begin{proof}
As $C(X)\rtimes \Gamma$ is finitely generated, we can apply \cite[Theorem 0.1]{MR1443393}.  It is therefore sufficient to construct semiprojective algebras $A_n$ and homomorphisms $\beta_n: C(X)\rtimes \Gamma\to A_n$ and $\beta'_n: A_j\to C(X)\rtimes \Gamma$ such that $\beta'_n\circ \beta_n(a)\to a$ with respect to the norm on $C(X)\rtimes \Gamma$.  We do this in much the same way as in the afforementioned theorem.  

To that end, find clopen covers $\mathcal{C}_n$ of $X$ made up of sets of diameter $<\frac{1}{n}$ such that $\alpha(\gamma)(C)\in \mathcal{C}_n$ for all $C\in \mathcal{C}_n$ and $\gamma\in \Gamma$.  Find finite subsets $X_n\subset X$ which are closed under the action by $\Gamma$ and have exactly one element in each $\mathcal{C}_n$.  

The inclusion maps $i_n: X_n\xhookrightarrow{} X$ induce dual maps $i_{n,*}: C(X)\to C(X_n)$, and since the $i_n$ are equivariant, so are the $i_{n, *}$.  Moreover, we can construct $p_n: X\to X_n$ by defining $p_n(C)$ to be the unique $x\in X_n$ such that $x\in C$.  These maps are similarly equivariant, as are their dual maps $p_{n, *}: C(X_n)\to C(X)$.  

By universality, the maps $i_{n, *}$ and $p_{n, *}$ induce homomorphisms $\beta_n: C(X)\rtimes \Gamma \to C(X_n)\rtimes \Gamma$ and $\beta'_n: C(X_n)\rtimes \Gamma\to C(X)\rtimes \Gamma$.  Since the algebras $C(X_n)\rtimes \Gamma$ are finite dimensional, they are semiprojective \cite[Corollary 2.28]{Blackadar1985}; and a continuity argument shows that $p_{n, *}(i_{n, *}(f))\to f$ uniformly so that $\beta'_n\circ \beta_n(a)\to a$ for all $a\in C(X)\rtimes \Gamma$.  \end{proof}

\example{Suppose $\Gamma$ is a finite group and $\Gamma\acts S^1$ where $S^1$ is the unit circle.  Then $C(S^1)\rtimes \Gamma$ is a $1$-NCCW and hence semiprojective.  }

\begin{proof}
Observe that $C(S^1)\rtimes \Gamma$ fits into the pullback diagram below:
\begin{center}
\begin{tikzcd}
  C(S^1)\rtimes \Gamma \arrow[r, "\phi"] \arrow[d, "\eta"]
    & C^*\Gamma \arrow[d, "\psi"] \\
  C([0, 1], C^*\Gamma) \arrow[r, "\text{ev}_0\oplus \text{ev}_1"]
&C^*\Gamma\oplus C^*\Gamma \end{tikzcd}
\end{center}

\noindent where $\eta(fu_\gamma)$ is the function $x\mapsto f(x)u_\gamma$, $\phi(fu_\gamma) = f(0)u_\gamma$, and $\psi$ is the diagonal map.  The algebra $C(S^1)\rtimes \Gamma$ is therefore a $1$-dimensional non-commutative CW complex in the sense of \cite{+1998+101+143}, and therefore semiprojective by \cite[Theorem 6.2.2]{+1998+101+143}.  \end{proof}

\subsection{Structure theory of virtually free groups}

\normalfont Section 4 will make use of the structure theory of virtually-free groups in much the same way as \cite{eilers}.  Specifically, we need a result from \cite{HERFORT1999511}, which we restate below.  

\theorem{If $\Gamma$ is a finitely generated, virtually-free group, then there is a group $G$ which is a finite tree product of finite groups such that $\Gamma$ is an iterated HNN extension of $G$ where each extension is taken over a finite subgroup.  }\label{virtually free structure theory}

\normalfont For the definition of a tree product, we refer to \cite{KIM2002323}.  For HNN extensions, one can see, for instance, \cite{eilers}.

\section{Almost-actions and perturbing problems}

\normalfont We begin with a definition based on our motivating examples, which we introduce in \ref{motivating examples lemma}. 

\definition{Let $\Gamma$ be a group and $X$ a compact metric space.  For $n\in \N$, let $E_n\subset X$ be a finite subset and $\alpha_n:\Gamma\to \text{Sym}E_n$.  We say the sequence $\mathbold{\alpha}= (\alpha_n)_{n\in \N}$ is an almost-action if 

$$\max_{e\in E_n} d_X(\alpha_n(\gamma)\circ \alpha_n(\delta)(e), \alpha_n(\gamma\delta)(e))\to 0$$
as $n\to \infty$ for all $\gamma, \delta\in \Gamma$.  We say $\mathbold{\alpha}$ is a uniform almost-action if 

$$\sup_{\gamma, \delta\in \Gamma}\max_{e\in E_n} d_X(\alpha_n(\gamma)\circ \alpha_n(\delta)(e), \alpha_n(\gamma\delta)(e))\to 0.$$  

\noindent If $S\subset \Gamma$ is a finite generating set, we say $\mathbold{\alpha}$ is a semi-uniform almost-action if

$$\sup_{\gamma\in \Gamma, s\in S} \max_{e\in E_n} d_X(\alpha_n(s)\circ \alpha_n(\gamma)(e), \alpha_n(s\gamma)(e))\to 0.$$

Notice that an almost-action of a finite group is a uniform almost-action.  Also, we can always assume the $\alpha_n$ respect inverses and the identity (just replace $\alpha_n(\gamma^{-1})$ with $\alpha_n(\gamma)^{-1}$).  If $\mathbold{\alpha}$ is an almost-action such that there is an action $\Gamma\acts X$ such that $\max_{e\in E_n}d_X(\gamma\cdot e, \alpha_n(\gamma)(e))\to 0$ for all $\gamma\in \Gamma$, we say $\mathbold{\alpha}$ is an almost-action approximating $\Gamma\acts X$, and similarly say it is an approximating almost-action uniformly approximating $\Gamma\acts X$ if the estimate holds independent of $\gamma$.  

If there exist (actual) actions $\tilde{\alpha}_n:\Gamma\to \text{Sym}E_n$ such that 

$$\max_{e\in E_n}d_X(\alpha_n(\gamma)(e), \tilde{\alpha}_n(\gamma)(e))\to 0$$
as $n\to\infty$ for all $\gamma\in \Gamma$ (respectively, $\sup_{\gamma\in \Gamma}\max_{e\in E_n}d_X(\alpha_n(\gamma)(e), \tilde{\alpha}_n(\gamma)(e))\to 0$), then we say the dynamical perturbing problem associated to $\mathbold{\alpha}$ has a solution (respectively, a uniform solution).  

If $\alpha_n: \Gamma\to \text{Sym}E_n$ and $\alpha_n': \Gamma\to \text{Sym}E'_n$ are two almost-actions such that there exist bijections $\phi_n: E_n\to E'_n$ such that $\sup_{e\in E_n} d_X(\phi_n(e), e)\to 0$ as $n\to\infty$ and such that $\sup_{e\in E_n} d_X(\phi_n(\alpha_n(\gamma)(e)), \alpha'_n(\gamma)(\phi_n(e)))\to 0$ for all $\gamma\in \Gamma$, we say $\mathbold{\alpha}$ and $\mathbold{\alpha'}$ are perturbations of each other.  Notice that in this case, the dynamical perturbing problem associated to $\mathbold{\alpha}$ has a solution if and only if the perturbing problem associated to $\mathbold{\alpha'}$ has a solution.  }\label{main definition}

\normalfont
We conclude this section by producing some motivating examples.  

\lemma{Suppose $\Gamma$ is a finitely-generated group, and $X$ is a compact metric space with no isolated points.  Given an action $\Gamma\acts X$ which fixes a measure of full support, there exists an almost-action approximating $\Gamma\acts X$.  }\label{motivating examples lemma}
\begin{proof}
Let $S$ be a finite generating set for $\Gamma$.  Let $F_S$ be the free group on the set $S$.  We can think of the action $\Gamma\acts X$ as an action $F_S\acts X$.  Apply \cite[Lemma 5.1]{kerr2011residually} to this action to get actions of $F_S$ on finite subsets $E_n\subset X$ which approximate $F_S\acts X$.  These actions are not in general actions of $\Gamma$, but they are almost-actions.  We can also prescribe that $E_n$ is $\epsilon$-dense in $X$ for any $\epsilon>0$ and $n$ sufficiently large.  \end{proof}

\lemma{Suppose $\alpha$ is an action of a group $\Gamma$ on a compact metric space $X$ and that $\alpha_n:\Gamma\to \text{Sym}(E_n)$ ($E_n\subset X$) is an almost-action which uniformly-approximates $\alpha$. Then $(\alpha_n)$ is a semi-uniform almost-action. }
\begin{proof}
For a given $\gamma\in \Gamma$, continuity of $\alpha(\gamma^{-1})$ implies that for any $\delta>0$ one can find $N$ such that

$$\alpha_n(\gamma^{-1}\gamma_0)(e) \approx_\delta \alpha(\gamma^{-1}\gamma_0)(e) = \alpha(\gamma^{-1})\circ \alpha(\gamma_0)(e) \approx_{\delta} \alpha(\gamma^{-1})\circ \alpha_n(\gamma_0)(e).$$

\end{proof}

\noindent \normalfont For approximating almost-actions coming from certain types of groups, we can improve the approximation.  It remains to be seen if further study of shadowing or other conditions may yield further refinements.  

\definition{Suppose $\Gamma$ is a finitely presented group.  We say $\Gamma$ has uniformly bounded combinatorial complexity if there exists a finite presentation $\langle S | R\rangle$, a set of normal forms $\{w_\gamma \mid \gamma\in \Gamma\}$ such that, for any $s\in S$, $\gamma\in \Gamma$, $sw_\gamma$ can be freely reduced to $w_{s\gamma}$ after multiplying on the left by one relation from the set $R\cup R^{-1}$.  }

\example{Suppose $\Gamma$ is finitely generated and virtually free.  Then $\Gamma$ has uniformly bounded combinatorial complexity.  }
\begin{proof}
We can assume $\Gamma$ has a free, normal subgroup of finite index.  So we have a short exact sequence $1\to F_n\to \Gamma\to G\to 1$ where $F_n$ is the free group on $n$ generators and $G$ is finite.  Let $T = \{g_1, \ldots, g_d\}$ be a set of coset representatives of $\Gamma/F_n$ (so $d = |G|$).  Let $f_1, \ldots, f_n$ be generators for $F_n$, and put $S = \{f_1, \ldots, f_n, g_1, \ldots, g_d\}$.  Then for each $\gamma\in \Gamma$, $\gamma$ can be written uniquely as $g_iw$ where $g_i\in T$ and $w$ is a freely reduced word in the $f_i$.  This will be our normal form for $\gamma$.  Let $\beta_i(f_j) = g_i^{-1}f_jg_i\in F_n$.  Further note that if the coset representative corresponding to $g_jg_i$ is $g_{(i, j)}$, then $g_{k}g_{i} = g_{(i, k)}f_{(i, k)}$ for some $f_{(i, j)}\in F_n$).  Let $R = \{\beta_i(f_j) = g_i^{-1}f_jg_i \mid 1\leq i\leq d\text{, }1\leq j\leq n\}\cup \{g_{k}g_{i} = g_{(i, k)}f_{(i, k)} \mid 1\leq i\leq d\text{, } 1\leq k\leq d\}$ (where we write $\beta_i(f_j)$ and $f_{i, k}$ as freely reduced words).  With the presentation $\langle S | R\rangle$, we now have the desired property.  \end{proof}

\lemma{Suppose $\Gamma$ has uniformly bounded combinatorial complexity with presentation $\langle S | R\rangle$ a presentation and $\{w_\gamma\}$ normal forms witnessing this property, and that $\mathbold{\alpha}$ is an almost-action of $\Gamma$.  Then there is a perturbation $\mathbold{\beta}$ of $\mathbold{\alpha}$ which is semi-uniform, i.e. $\sup_{\gamma\in \Gamma} d_X(\beta_n(s)\circ \beta_n(\gamma), \beta_n(f\gamma))\to 0$ for all $s\in S$.  }

\begin{proof}
For $w_\gamma = s_k\cdots s_1$, define $\beta_n(\gamma) = \alpha_n(w_\gamma) := \alpha_n(s_k)\circ \cdots \circ \alpha_n(s_1)$.  We can assume $\alpha_n$ respects inverses, so free reduction is no problem.  Then find $N$ such that if $r\in R$ (recalling that $r$ is a word in $S$) we have $d_X(\alpha_n(r)(e), e)<\epsilon$ for $n\geq N$.  Then for $s\in S$, $\gamma\in \Gamma$, there is $r\in R$ such that the word $sw_\gamma$ is equal to the word $rw_{s\gamma}$, after some free reduction.  Hence,

$$\beta_n(s)\circ \beta_n(\gamma) = \alpha_n(s)\circ \alpha_n(w_\gamma) \approx_{\epsilon} \alpha_n(r)\circ \alpha_n(w_{s\gamma}) = \beta_n(s\gamma)$$ \end{proof}

\section{Elementary solutions when X is a Cantor set}

\normalfont
\noindent
When $X$ is a Cantor space, we can find solutions to many dynamical perturbing problems by elementary means.  

\lemma{Let $X$ be the Cantor space and let $\Gamma = \Gamma_0*_\Delta \Lambda$ be a group with $\Lambda$ a finite group and $\Delta$ subgroup of $\Lambda$.  Suppose $\Gamma\acts^{\beta} X$ and $\mathbold{\alpha} =(\alpha_n)$ an almost-action approximating $\beta$.  If the dynamical perturbing problem associated to $(\alpha_{n|\Gamma_0})$ has a solution, then the dynamical perturbing problem associated to $\mathbold{\alpha}$ has a solution.  }
\begin{proof}

Consider the approximating almost-action $(\alpha_{n|\Lambda})$.  Using our solution to the perturbing problem associated to $(\alpha_{n|\Gamma_0})$, we can replace $\alpha_{n|\Lambda}(\delta)$ by the action of $\delta$ according to the solution.  This yields a perturbation of the original almost-action, so our new almost-action $(\alpha_{n|\Lambda})$ will still approximate $\beta$ and will now be an action when restricted to $\Delta$.  We are therefore done if we can find a solution to the perturbing problem associated to $(\alpha_{n|\Lambda})$ which doesn't change $(\alpha_{n|\Delta})$

Fix $\epsilon>0$ and find a clopen partition $\mathcal{X}$ of $X$ which is preserved by $\beta(\Lambda)$ and whose elements have diameter $<\epsilon$.  Let $\delta>0$ be the least distance between two elements of $\mathcal{X}$, and find $N$ such that for $n\geq N$, $d_X(\alpha_n(\lambda)(e), \beta(\lambda)(e))<\delta$ for all $\lambda\in \Lambda$.  Then whenever $\beta(\lambda)$ maps $X_i\in \mathcal{X}$ to $X_j\in \mathcal{X}$, and $e\in E_n\cap X_i$, $\alpha_n(\lambda)(e)\in X_j$.  Hence, $\# E_n\cap X_i = \# E_n\cap X_j$.  

Choose a labelling $1, \ldots, \#E_n\cap X_i$ for the elements of each $E_n\cap X_i$.  It is possible to do this so that the labellings of $e$ and $\beta(\delta)(e)$ are the same for each $\delta\in \Delta$ and $e\in E_n$.  Now for $e\in X_i$ define $\tilde{\alpha}_n(\gamma)(e)$ to be the unique $e'\in \beta(\gamma)(X_i)$ with the same label as $e$, and notice that this defines an action and doesn't change $\alpha_{n|\Lambda}$.  By the previous paragraph, this is possible, and changes $\alpha_n$ by less than $\epsilon$.  \end{proof}

%

\lemma{Suppose $\Gamma$ is a finite tree product of finite groups, and $\Gamma\acts^\alpha X$.  If $\mathbold{\alpha} = (\alpha_n)$ is an almost-action approximating $\alpha$, the dynamical perturbing problem associated to $\mathbold{\alpha}$ has a solution.  }

\begin{proof}
A finite tree product of finite groups is constructed by taking iterated amalgamated free products over finite groups, starting with a finite group.  The lemma therefore follows from the previous lemma and induction.  \end{proof}

\lemma{Let $X$ be the Cantor space and let $\Gamma\acts X$ be an action by homeomorphisms.  If $\Gamma = *_\phi \Lambda$ is an HNN extension of a group $\Lambda$ over a finite subgroup (meaning $\phi: F\to G$ is an isomorphism of finite subgroups).  Suppose $\mathbold{\alpha}$ is an approximating almost-action such that the perturbing problem associated to $\mathbold{\alpha}_{|\Lambda}$ is solvable.  Then the perturbing problem associated to $\mathbold{\alpha}$ is solvable.  }

\begin{proof}
Take a solution $\widetilde{\mathbold{\alpha}_{|\Lambda}}$ to the perturbing problem associated to $\mathbold{\alpha}_{|\Lambda}$, and replace $\alpha_n(\lambda)$ by $\widetilde{\alpha_n}(\lambda)$ for all $\lambda\in \Lambda$.  This yields a perturbation of $\mathbold{\alpha}$.  

Fix $\epsilon>0$ and find clopen partitions $\mathcal{X}$ and $\mathcal{Y}$ of $X$ preserved by $F$ and $G$ respectively whose elements have diameter at most $\epsilon$.  Let 

$$\delta =\min \{\min_{U, V\in \mathcal{X}} d_X(U, V), \min_{U, V\in \mathcal{Y}} d_X(U, V)\}.$$
Notice that $h$ takes $U\in \mathcal{X}$ to $V\in \mathcal{Y}$, so if $\alpha_n(h)$ approximates the action of $h$ to within $<\delta$, $\alpha_n(h)(e)\in V$ for any $e\in U$.  This implies $\#E_n\cap U = \#E_n\cap V$.  We can choose a labelling $1, \ldots, \#E_n\cap U$ for each $U\in \mathcal{X}$ in such a way that the labellings are preserved by the action of $F$, and another labelling with the same properties for each $V\in \mathcal{Y}$ and the action of $G$.  Then for $e\in U$, define $\tilde{\alpha}_n(h)(e)$ to be the $e'\in h\cdot U$ with the same label as $e$ (this is possible since $\#E_n\cap U = \#E_n\cap V$, and changes what $\alpha_n(h)$ does by at most $\epsilon$).  We now have an action $\Gamma\acts E_n$.  \end{proof}

\theorem{Suppose $X$ is a Cantor set, $\Gamma$ is finitely generated and virtually free, and $\mathbold{\alpha}$ is an almost-action which approximates an action $\Gamma\acts X$.  Then the dynamical perturbing problem associated to $\mathbold{\alpha}$ has a solution.  }\label{cantor sets theorem}

\begin{proof}
This follows from the structural result \ref{virtually free structure theory}, as every virtually free group is obtained by beginning with a finite tree product of finite groups (meaning taking finitely-many amalgamated products over finite groups) and then taking iterated HNN extensions over finite groups.  \end{proof}

\corollary{Suppose $X$ is a Cantor set, $\Gamma$ is virtually free, and $\Gamma\acts X$ is an action by homeomorphisms which fixes a measure of full support.  Then $\Gamma\acts X$ is residually finite.  }\label{main corollary}
\begin{proof}
By \ref{motivating examples lemma}, there exists an approximating almost-action approximating $\Gamma\acts X$, and so the result follows from the previous theorem.  \end{proof}

\section{A representation}

\normalfont In light of earlier results about asymptotic representations and Gelfand duality, it is natural to consider whether approximate actions can be used to construct asymptotic representations.  Going about this in a na\"{i}ve way, however, almost-actions as defined here lead to operators which only give an asymptotic representation in a pointwise sense, since the error is dependent on the continuity of a function.  One potential way around this is to focus on the case of finite groups, but even then it's not clear whether the asymptotic representations will be relevant to their associated dynamical perturbing problems (for instance, conjugation by the lifted unitaries won't necessarily preserve the diagonal).  Another potential way forward is to look for asymptotic representations of the crossed product $C^*$-algebra as opposed to the group itself, so that each element has a fixed element of $C(X)$ attached to it.  

In what follows, the left regular representation of $\Gamma$ will be denoted $\lambda$, and, for $f\in C(X)$, the representation of $f$ as a multiplication operator on $L^2(X)$ will be denoted $M_f$.


%
%
%
%

\lemma{Suppose $\alpha: \Gamma\to \text{Homeo}(X)$ is an action and $\widehat{\alpha}$ is the dual action on $C(X)$ given by $\widehat{\alpha}(\gamma)(f)(x) = f(\alpha(\gamma^{-1})(x))$.  Suppose $\mathbold{\alpha} = (\alpha_n: \Gamma\to \text{Sym}E_n)$ is a semi-uniform almost-action approximating $\alpha$ (for instance one associated to $\alpha$ through \cite[Lemma 5.1]{kerr2011residually} when $\Gamma$ is virtually free).  Define $[\widehat{\alpha}_n(\gamma)(f)](e) = f(\alpha_n(\gamma^{-1})(e))$, that is, the dual almost-action induced from $\alpha_n$.  

There is a representation $\rho_\mathbold{\alpha}:C(X)\rtimes \Gamma\to \frac{\prod_n \mathcal{B}(l^2E_n \otimes l^2\Gamma)}{\bigoplus_n \mathcal{B}(l^2E_n\otimes l^2\Gamma)}$ coming from $(1\otimes \lambda)$ (where $\lambda$ is the regular representation of $\Gamma$) and representations $\pi_n: C(X)\to B(l^2E_n\otimes l^2\Gamma)$ defined by $\pi_n(f)(\xi\otimes \delta_\gamma) = M_{\widehat{\alpha}_n(\gamma^{-1})(f|_{E_n})}(\xi)\otimes \delta_\gamma$.  } \label{representation}
\begin{proof}

The definition of $\widehat{\alpha}_n$ gives set maps $\widehat{\alpha}_n:\Gamma\to \text{Aut}(C(E_n))$.  Let $\pi_n$ be defined as above.  Notice that $(1\otimes \lambda_\gamma)(\pi_n(f))(1\otimes \lambda^*_{\gamma})$ takes $\xi\otimes \delta_{\gamma_0}$ to $M_{\widehat{\alpha}_n(\gamma_0^{-1}\gamma)(f|_{E_n})}(\xi)\otimes \delta_{\gamma_0}$ whereas $\pi_n(\widehat{\alpha}(\gamma)(f))$ takes $\xi\otimes \delta_{\gamma_0}$ to $M_{\widehat{\alpha}_n(\gamma_0^{-1})(\widehat{\alpha}(\gamma)(f)|_{E_n})}(\xi)\otimes \delta_{\gamma_0}$.  

Fix $\epsilon>0$.  For a given $\gamma\in \Gamma$ and any $\delta>0$ one can find $N$ such that 
$$\alpha_n(\gamma^{-1}\gamma_0)(e) \approx_\delta \alpha_n(\gamma^{-1})\circ \alpha_n(\gamma_0)(e) \approx_\delta \alpha(\gamma)^{-1}\circ \alpha_n(\gamma_0)(e)$$


\noindent for any $e\in E_n$, $\gamma_0\in \Gamma$, and all $n\geq N$ (using that our almost-action is semi-uniform).  Then for a given $f\in C(X)$ one can let $\delta$ be small enough that continuity of $f$ implies

$$\|\widehat{\alpha}_n(\gamma_0^{-1}\gamma)(f|_{E_n}) - \widehat{\alpha}_n(\gamma_0^{-1})(\widehat{\alpha}(\gamma)(f)|_{E_n})\|_{C(E_n)} < \epsilon$$

\noindent for sufficiently large $n$, and so

$$\|(1\otimes \lambda_\gamma)(\pi_n(f))(1\otimes \lambda^*_{\gamma}) - \pi_n(\widehat{\alpha}(\gamma)(f))\|_{\mathcal{B}(l^2E_n\otimes l^2\Gamma)}<\epsilon\text{ .}$$

In other words, the failure of covariance between $(1\otimes \lambda_\gamma)$ and $\pi_n(f)$ tends to zero as $n\to\infty$.  We can therefore build a representation of $C(X)\rtimes\Gamma$ into $\frac{\prod_n \mathcal{B}(l^2E_n \otimes l^2\Gamma)}{\bigoplus_n \mathcal{B}(l^2E_n\otimes l^2\Gamma)}$.  \end{proof}

\remark{Notice that the restriction of $\rho_{\mathbold{\alpha}}$ to either $C^*\Gamma$ or $C(X)$ is already induced by a representation into $\prod_n\mathcal{B}(l^2E_n\otimes l^2\Gamma)$.  }



\normalfont

We are more interested in the perturbing problems arising from non-uniform almost-actions, as these are the almost-actions of which we have naturally-arising examples.  It should be noted however that, even for finite groups, it is not obvious how to solve these perturbing problems through elementary means when $X$ is not a Cantor set.  Even producing solutions for finite groups would immediately give solutions for free products of finite groups and possibly even amalgamated products.

\section{Almost-actions and conditional lifting}

\normalfont We can use some of the basic background from the theory of Bratteli diagrams \cite[Theorem 2.3.5]{brattelidiagrams} to show that morphisms of finite dimensional $C^*$-algebras are conditionally semiprojective.  First let us define what that means. 

\definition{A homomorphism of $C^*$-algebras $\phi: A\to B$ is conditionally semiprojective if, whenever we have a representation $\rho: B\to C/\overline{\cup_n J_n}$ (for $\cup_n J_n$ an increasing union of ideals) and a lift of $\rho\circ \phi$, $\widetilde{\rho\circ \phi}: A\to C/J_m$ for some $m$, there is a lift $\tilde{\rho}: B\to C/J_{m'}$ ($m'\geq m$) of $\rho$ such that $\tilde{\rho}\circ \phi =\pi_{m'} \circ \widetilde{\rho\circ \phi}$.  Conditional weak semiprojectivity and conditional matricial weak semiprojectivity are defined similarly.}\label{conditional semiprojectivity def}

\lemma{Suppose $\rho: \mathcal{A}\to \mathcal{B}/\overline{\cup_n J_n}$ is a representation.  One can therefore choose an asymptotic representation $\rho_n: \mathcal{A}\to \mathcal{B}/J_n$.  Suppose $v_1, \ldots, v_k\in \mathcal{A}$ are orthogonal partial isometries such that $\rho_n(v_1) + \cdots + \rho_n(v_k) = v$, a partial isometry in $\mathcal{B}/J_n$ for all $n$.  Let $\epsilon>0$.  Then for all $n$ sufficiently large, there are $\tilde{v}_1, \ldots, \tilde{v}_k\in \mathcal{B}/J_n$ orthogonal partial isometries such that $\|\tilde{v}_i - \rho_n(v_i)\|<\epsilon$ and $\tilde{v}_1 + \cdots \tilde{v}_k = v$.  }\label{orthogonal lifting}

\begin{proof}
Apply \ref{lifting partial isometries} to each $\rho_n(v_i)$ to obtain orthogonal partial isometries $\bar{v}_i$ in each $\mathcal{B}/J_n$.  Then for $n$ sufficiently large, $\bar{v}_1 + \cdots + \bar{v}_k \approx_{\epsilon} v$.  Since these two partial isometries are close, their source and range projections must be close.  Using the same trick as in \ref{lifting partial isometries} to find unitaries which conjugate one projection over to the other, we have unitaries $u_1$ and $u_2$ which are close to the identity such that $u_2(\bar{v}_1 + \cdots + \bar{v}_k)u_1 = v$.  Setting $\tilde{v}_i = u_2\bar{v}_iu_1$ completes the proof.  \end{proof}

\lemma{A canonical homomorphism of finite dimensional $C^*$-algebras is conditionally semiprojective.  }

\begin{proof}

Let $h: F_1\to F_2$ be such a map and suppose $\rho_2: F_2\to \mathcal{B}/\overline{\cup J_n}$ and $\rho_1: F_1\to \mathcal{B}/J_n$ are representations such that $\pi \circ \rho_1 = \rho_2$.  Assume $F_2\cong \bigoplus_{n=1}^N M_{k_n}(\C)$ and fix matrix units $e_{ij}^{(n)}$ for each $M_{k_n}(\C)$ ($1\leq i, j\leq k_n$).  Our goal is to define a lift $\tilde{\rho}_2$ by replacing the image of each matrix unit with a partial isometry so that the appropriate relations are satisfied \textit{and} so that the restriction of $\tilde{\rho}_2$ to the $h(F_1)$ is still the same as $\rho_1$.  

Since $h$ is canonical, the image of a matrix unit in $F_1$ under $h$ is a sum of orthogonal matrix units in $F_2$ of the form $e^{(n)}_{i,j} + e^{(n)}_{i+l, j+l} + \cdots + e^{(n)}_{i+ml, j+ml}$ for some integers $l$ and $m$.  The image of each of these matrix units under $\rho_2$ is close to being a partial isometry.  By \ref{orthogonal lifting}, we can lift the image under $\rho_2$ of each matrix unit in this sum while keeping the sum of their images the same.  We can then choose lifts of other matrix units in $F_2$ until those lifts determine lifts for all matrix units of $F_2$.  \end{proof}

\theorem{Suppose $\phi: A\to B$ is a homomorphism of finite-dimensional $C^*$-algebras.  Then $\phi$ is conditionally semiprojective.  }\label{finite dimensional algebras cond semiprojective}
\begin{proof}
By possibly replacing $A$ by the quotient $A/\ker(\phi)$, we can assume $\phi$ is injective.  Every injective morphism of finite-dimensional $C^*$-algebras is inner-equivalent to a canonical homomorphism, that is, there are inner automorphisms $h_1: A\to A$ and $h_2: B\to B$.  We therefore have the following commutative diagram   

\begin{center}

\begin{tikzcd}
  A \arrow[r, "\phi"] \arrow[d, "h_1"]
    & B \\
  A \arrow[r, "\psi"]
&B \arrow[u, "h_2^{-1}"]  \end{tikzcd}

\end{center}

\noindent where $\psi$ is a canonical homomorphism.  By the previous lemma, $\psi$ is conditionally semiprojective, from which it follows that $\phi$ is.  \end{proof}

\normalfont \noindent We can also see that asking for conditional lifts of the asymptotic representations constructed in the previous section is too much.  

\example{The inclusion $C(S^1)\hookrightarrow C(S^1)\rtimes F$ for any (non-trivial) action of a finite group $F$ on $S^1$ is not conditionally matricially weakly semiprojective.  In particular, this gives an example of an inclusion of a sub-$1$-NCCW (as a Cartan subalgebra) into a $1$-NCCW which is not conditionally matricially weakly semiprojective.  }\label{counter-example}

\begin{proof}
Let $(\alpha_n)$ be an approximating almost-action on $E_n\subset X$ constructed for instance in \ref{motivating examples lemma} (in fact, since $\Gamma$ is finite, constructing such an almost-action is trivial, and we could even have it be an honest action by simply perturbing an orbit of $\Gamma\acts X$).  There are many possible choices for what the sets $E_n$ can be.  In particular, any perturbation of $(\alpha_n)$ will be another approximating almost-action.  Let $D_n = \rho_n(C(S^1))$.  If the inclusion $C(S^1)\to C(S^1)\rtimes F$ is conditionally matricially weakly semiprojective, then the representation $\rho$ coming from \ref{representation} has a lift $\tilde{\rho}$ such that $\tilde{\rho}_n(C(S^1)) = D_n$ for $n$ sufficiently large.  Then we would have an action of $\Gamma$ on $D_n$ (via conjugation by $\tilde{\rho}_n(u_\gamma)$).  Now, $D_n$ is isomorphic to $\C^{E_n}$ and moreover the composition of $\tilde{\rho}_n$ restricted to $C(S^1)$ with this isomorphism is simply the evaluation map $C(S^1)\to \C^{E_n}$.  Pushing the action $\Gamma\acts D_n$ forward to $\C^{E_n}$ via this isomorphism, we get an action $\Gamma\acts \C^{E_n}$ such that the evaluation map $C(S^1)\to \C^{E_n}$ is equivariant.  Then by duality, we would have that the inclusion $E_n\to S^1$ is equivariant.  But that means the orbit under $\Gamma\acts S^1$ of $x\in E_n$ must be contained in $E_n$.  By perturbing all but one point of $E_n$ (and using that the action is non-trivial), we can therefore produce many examples of representations $\rho$ where such a lift cannot exist.  \end{proof}

\section{Lifting Cartan structures}\label{lifting Cartan pairs}

\normalfont We have seen that asking for a lift of the representation $\rho_{\mathbold{\alpha}}$ constructed from an almost-action which doesn't change the natural lift of ${\rho_{\mathbold{\alpha}}}_{|C(X)}$ is too much, in that it leads to an equivariant map of the underlying dynamics.  It is also worth noting that the existence of a (non-conditional) lift of $\rho_{\mathbold{\alpha}}$ doesn't appear to give a solution to the perturbing problem associated to $\mathbold{\alpha}$.  In the case where $\Gamma$ is finite, such a lift can be seen to always exist, and this simply corresponds to restricting the original action to the union of the (actual) orbits of each $e\in E_n$.  

Notice that an approximately equivariant map $Y\to X$ induces an approximate homomorphism of crossed products $C(X)\rtimes \Gamma\to C(Y)\rtimes \Gamma$ with some extra properties: the restriction to either $C(X)$ or $C^*(\Gamma)$ is a homomorphism, and the image of the Cartan pair $(C(X)\rtimes \Gamma, C(X))$ is another Cartan pair.  The approximate representations constructed in \ref{representation}, have the former property, but the latter property only holds approximately.  To characterize the perturbing problem associated to $\mathbold{\alpha}$, what we actually want are homomorphisms $\phi_n: C^*\Gamma\to \mathcal{B}(l^2E_n \otimes l^2\Gamma)$ which are approximately covariant with $\pi_n$ and are such that $\phi_n(u_\gamma)$ normalizes $\pi_n(C(X))$.  Although this is enough by itself to characterize perturbing problems, the natural $C^*$-algebraic structure to consider in this situation is that of Cartan pairs, and we will see we can include this extra structure and still get a $C^*$-algebraic characterization in terms of trying to ``lift" the structure of a Cartan pair, given $\rho_{\mathbold{\alpha}}$ as constructed in \ref{representation}.  We will give a general definition of what it means to lift the structure of a Cartan pair, and bring in the full structure of said in the final section to characterize abstract perturbing problems involving finite groups.

\theorem{Suppose $\mathbold{\alpha}$ is an approximating almost-action such that we can construct the representation $\rho_{\mathbold{\alpha}}$ (e.g. if $\mathbold{\alpha}$ is semi-uniform or uniformly-approximating).  Let $\rho_\mathbold{\alpha}: C(X)\rtimes \Gamma\to \frac{\prod_n \mathcal{B}(l^2E_n \otimes l^2\Gamma)}{\bigoplus_n \mathcal{B}(l^2E_n\otimes l^2\Gamma)}$ be the representation constructed in \ref{representation} and $(\pi_n)$ the canonical lift of $\rho_\mathbold{\alpha}$ restricted to $C(X)$.

The dynamical perturbing problem associated to $\mathbold{\alpha}$ has a solution if and only if there is a homomorphism $\phi:C^*\Gamma\to \prod_n \mathcal{B}(l^2E_n \otimes l^2\Gamma)$ such that 

\begin{description}
\item(i) $(\pi_n)$ and $(\phi_n)$ are approximately-covariant
\item (ii) $\Big(\phi_n(C^*\Gamma)\cup \pi_n(C(X)), \pi_n(C(X))\Big)$ is a Cartan pair for all $n$. 

\end{description} }\label{first Cartan characterization}

\begin{proof}

Suppose we have such a homomorphism.  There is a homomorphism of groups $\Gamma\to \Gamma'$ underlying $\phi$, so $\phi(C^*\Gamma)$ is minimally generated by a linearly independent subset $S\subset \phi(\{u_\gamma \mid \gamma\in \Gamma\})$ which generates $\Gamma'$ as a group.  Then by $(ii)$, $S$ can be chosen such that every element of $S$ normalizes $\psi_n(C(X))$, so $\phi_n(u_\gamma)$ normalizes $\psi_n(C(X))$ for any $\gamma$ and $n$.  Since $\psi_n(C(X)) = \pi_n(C(X))\cong \mathbb{C}^{E_n}$ for all $n$, we have, by Gelfand duality, actions $\Gamma\acts E_n$.  Since $\psi$ and $\phi$ give a covariant representation $C(X)\rtimes \Gamma\to \frac{\prod_n \mathcal{B}(l^2E_n \otimes l^2\Gamma)}{\bigoplus_n \mathcal{B}(l^2E_n\otimes l^2\Gamma)}$, we get approximate equivariance on the level of topology, so these new actions approximate the original action $\Gamma\acts X$ (and so are a perturbation of the almost-action $(\alpha_n)$).  They therefore give a solution to the dynamical perturbing problem associated to $\mathbold{\alpha}$.

Now suppose we have a solution, $(\tilde{\alpha}_n)$, to the perturbing problem associated to $\mathbold{\alpha}$.  Define $\phi_n(u_\gamma)$ to be the block-diagonal operator in $\mathcal{B}(l^2E_n\otimes l^2\Gamma)$ with the permutation matrix corresponding to the permutation $\alpha_n(\lambda)\circ \tilde{\alpha}_n(\gamma)\circ \alpha_n(\lambda)^{-1}$ in the block corresponding to $\lambda\in \Gamma$.  This defines a homomorphism on $C^*\Gamma$.  Then, the map $\prod_n\mathbb{C}^{E_n}\to \prod_n\mathcal{B}(l^2E_n\otimes l^2\Gamma)$ given by $(v_n)\mapsto (w_n\otimes \delta_\lambda\to M_{\alpha_n(\lambda)^{-1}(v_n)}w_n\otimes \delta_\lambda)$ is equivariant for the actions given by $(\tilde{\alpha}_n)$ and conjugation by $\phi_n(u_\gamma)$.  This implies that if $v\in \mathbb{C}^{E_n}$ with $v = \psi_n(f) = \pi_n(f)$ and if $g\in C(X)$ is such that $g(e) = \tilde{\alpha}_n(\gamma)(f)(e)$ (such a $g$ always exists), then $\phi_n(u_\gamma)^*\psi_n(f)\phi_n(u_\gamma) = \psi_n(g)$.  We therefore see that $C^*(\psi_n(C(X)), \phi_n(C^*\Gamma))$ is isomorphic to the algebra generated by $\mathbb{C}^{E_i}$ (thought of as diagonal matrices) and the set of matrices $\Sigma_n$ corresponding to the permutations $\{\tilde{\alpha}_n(\gamma) \mid \gamma\in \Gamma\}$, and so $(C^*(\psi_n(C(X)), \phi_n(C^*\Gamma)), \psi_n(C(X)))\cong (C^*(\C^{E_i}, \Sigma_n), \C^{E_i})$ is a Cartan pair.


Since $\phi_n(u_\gamma)$ and $\psi_n = \pi_n$ give an approximately-covariant representation $C(X)\rtimes\Gamma\to \mathcal{B}(l^2E_n\otimes l^2\Gamma)$, we get a homomorphism $C(X)\rtimes\Gamma\to \frac{\prod_n \mathcal{B}(l^2E_n \otimes l^2\Gamma)}{\bigoplus_n \mathcal{B}(l^2E_n \otimes l^2\Gamma)}$. \end{proof}

\section{Abstract perturbing problems}
\normalfont

\noindent The characterization given in \ref{lifting Cartan pairs} relies on the almost-action approximating some actual action.  In this section, we will give a similar characterization for much more general perturbing problems when the group $\Gamma$ is finite.  We do this by constructing a commutative $C^*$-algebra and $\Gamma$-action from the almost-action, using a similar idea as in the construction of the regular representation of an action.  Like in \ref{representation}, this construction essentially trades not being an action for not preserving a subalgebra.  Notice, for instance, that if $\mathbold{\alpha}$ were a sequence of genuine actions, that permuting the blocks of $\rho_{\mathbold{\alpha}}(f)$ according to $\gamma$ is the same as doing the permutation associated to $\gamma$ in each block.  We will use a similar idea of putting together the ``orbits" of the almost-action into a regular representation.  

The following language is not necessary to state our main result in this section, but allows us to give some exposition and context in terms of general (unital) Cartan pairs as opposed to restricting attention only to crossed products. 

\definition{Suppose $\mathcal{A}$ is a $C^*$-algebra.  A subalgebra $\mathcal{B}\subset \mathcal{A}$ is a Cartan subalgebra, and the pair $(\mathcal{A}, \mathcal{B})$ is a Cartan pair if the following hold:
\begin{description}
\item(i) $\mathcal{B}$ is maximal abelian in $\mathcal{A}$
\item(ii) the normalizer of $\mathcal{B}$ (i.e. the set $\{n\in \mathcal{A} \mid nbn^{-1}\in \mathcal{B}\text{ for all }b\in \mathcal{B} \}$) generates $\mathcal{A}$
\item(iii) there exists a faithful conditional expectation $E:\mathcal{A}\to \mathcal{B}$.  
\end{description}}

\definition{Let $(\mathcal{A}, \mathcal{B})$ be a unital Cartan pair with conditional expectation $E:\mathcal{A}\to \mathcal{B}$.  Let $N = \{n\in \mathcal{A} \mid nbn^{-1}\in \mathcal{B}\text{ for all }b\in \mathcal{B} \}$ be the normalizer of $B$ and define $\mathcal{N} = C^*(E^{-1}(\{\mathds{1}_{\mathcal{B}}\})\cap N)$.  We call this the normalizer algebra of $(\mathcal{A}, \mathcal{B})$.  Observe that $C^*(\mathcal{B}, \mathcal{N}) = \mathcal{A}$, and that representations of $\mathcal{B}$ and $\mathcal{N}$ into the same algebra $\mathcal{C}$ which are covariant in the obvious sense give rise to a representation $\mathcal{A}\to \mathcal{C}$.  Furthermore, the set $\{u_\gamma \mid \gamma\in \Gamma\}\subset C(X)\rtimes \Gamma$ are linearly independent as a $C(X)$-vector space, and so representations $\psi$ and $\phi$ of $C(X)$ and $C^*\Gamma$ always give rise to a $*$-linear map of $C(X)\rtimes \Gamma$.}

\example{If $\mathcal{A} = C(X)\rtimes \Gamma$ is the (full or reduced) crossed product of an action $\Gamma\acts X$ and $B = C(X)$, then $\mathcal{N} = C^*(\{u_\gamma \mid \gamma\in \Gamma\})$.  Aside from the maximal abelian condition (which is satisfied, for instance, if the action is free) $C(X)$ satisfies all the requirements of being a Cartan subalgebra of $C(X)\rtimes \Gamma$.  }


\normalfont The essence of our problem is the following: If we have a map $\mathcal{A}\to\mathcal{C}/I$ which takes a Cartan pair to a Cartan pair (so as a map $\mathcal{A}\to \mathcal{C}$, it is close to taking a Cartan pair to a Cartan pair in some sense), can we ``wiggle" that map (by a small amount in some sense) and get a map $\mathcal{A}\to \mathcal{C}$ with the same property?  From a slightly different, more narrow, point of view, if $(\mathcal{A}, \mathcal{B})$ is a Cartan pair with normalizer algebra $\mathcal{N}$, and we have homomorphisms $\psi_n:\mathcal{B}\to \mathcal{C}_n$ and $\phi_n:\mathcal{\mathcal{N}}\to \mathcal{C}_n$ such that $\psi_n(\mathcal{B})$ is close to being Cartan in the subalgebra $C^*(\psi_n(\mathcal{B}), \phi_n(\mathcal{N}))$, can we ``wiggle" these maps a little bit and get homomorphisms $\widetilde{\psi}_n$ and $\widetilde{\phi}_n$ such that $\widetilde{\psi}_n(\mathcal{B})$ is Cartan in $C^*(\psi_n(\mathcal{B}), \phi_n(\mathcal{N}))$.  A priori, there are many reasonable definitions one could write down to capture this idea. We will restrict attention to one, informed by our goal of characterizing dynamical perturbing problems.

\definition{Suppose $(\mathcal{A}, \mathcal{B})$ is a Cartan pair, that $\Phi:\mathcal{A}\to \mathcal{C}/\mathcal{I}$ is a homomorphism such that $\big(\Phi(\mathcal{A}\big), \Phi(\mathcal{B}))$ is a Cartan pair, and that $\psi$ is a lift of $\Phi_{|\mathcal{B}}$.  We say the Cartan lifting problem associated to $(\mathcal{A}, \mathcal{B}, \Phi, \psi)$ has a solution if there is a $*$-linear map $\widetilde{\Phi}: \mathcal{A}\to \mathcal{C}$ such that

\begin{description}
\item(1) $\widetilde{\Phi}_{|\mathcal{B}} = \psi$
\item(2) $\pi \circ \widetilde{\Phi}: \mathcal{A}\to\mathcal{C}/\mathcal{I}$ is a homomorphism
\item(3) $\big(\widetilde{\Phi}(\mathcal{A}), \widetilde{\Phi}(\mathcal{B})\big)$ is a Cartan pair.

\end{description}}\label{Cartan lift definition}

\normalfont Classically, a ``lifting problem" should be something like the following.  If $A$ and $B$ are objects, $q:B \to B/I$ is a quotient, and $f$ is a map into $A\to B/I$ which preserves some structure, is there a map $\tilde{f}: A\to B$ which preserves the same structure and has the property that $q\circ \tilde{f} = f$.  We appear, then, to have departed somewhat from this spirit in the previous definition because we have not required $\widetilde{\Phi}$ to be a lift of $\Phi$; so, a priori, we have weakened the condition that $q\circ \tilde{f} = f$.  We will see in \ref{alt definition remark} that this is not really the case.

\lemma{Suppose $(\alpha_n)$ is an almost-action of a finite group $\Gamma$.  Then there is another almost-action $(\beta_n)$ such that $d_{X_n}(\beta_n(\gamma)(x), \beta(\delta)(x))>1$ for all $x\in X_n$ and $\gamma, \delta\in \Gamma$ with $\gamma\neq \delta$; and the dynamical perturbing problem associated to $(\beta_n)$ has a solution if and only if the one associated to $(\alpha_n)$ has a solution.  } \label{freeness lemma}

\begin{proof}
We will refer to elements in the disjoint union $\bigsqcup_{\gamma\in \Gamma} X_n$ by $x_\gamma$ for $x\in X_n$ and $\gamma\in \Gamma$.  Define $\beta_n: \Gamma\to \text{Homeo}(\bigsqcup_{\gamma\in \Gamma} X_n)$ by $\beta_n(\gamma)(x_\delta) = \alpha_n(\gamma)(x)_{\gamma\lambda}$.  Then, after choosing an appropriate metric for the disjoint union, $(\beta_n)$ has the required property.  It is an almost-action since $(\alpha_n)$ is.  

Assume the perturbing problem associated to $(\beta_n)$ has a solution $(\widetilde{\beta}_n)$.  Let $p: \bigsqcup_{\gamma\in \Gamma} X_n$ be the map given by $x_\gamma\mapsto x$.  Then define $\alpha_n(\gamma)(x) = p(\widetilde{\beta}(\gamma)(x_e))$.  This is an action because $(\widetilde{\beta}_n)$ is, and solves the perturbing problem associated to $(\alpha_n)$ because $(\widetilde{\beta}_n)$ solves the perturbing problem associated to $(\beta_n)$.  

Now assume the perturbing problem associated to $(\alpha_n)$ has a solution $(\widetilde{\alpha}_n)$.  Repeat the construction of $(\beta_n)$ with $(\widetilde{\alpha}_n)$ in place of $(\alpha_n)$.  This gives a solution to the perturbing problem associated to $(\beta_n)$\end{proof}


\theorem{Let $\Gamma$ be a finite group and $(X_n, d_{X_n})$ compact metric spaces such that $C(X_n)$ is generated for all $n$ by a set of at most $M$ elements, $\mathcal{F}_n$, and that $\bigcup_n \mathcal{F}_n$ is an equicontinuous family (e.g. each $X_n$ is a compact Riemannian manifold with dimension at most $d$ for some $d$).  

Let $\mathbold{\alpha} = (\alpha_n)$ with $\alpha_n: \Gamma\to \text{Homeo}(X_n)$ be an almost-action.  One can associate to $\mathbold{\alpha}$ a commutative $C^*$-algebra $\mathcal{B}_\mathbold{\alpha}$ and an action $\Gamma\acts \mathcal{B}_\mathbold{\alpha}$ with crossed product $\mathcal{A}_\mathbold{\alpha}$ such that 

\begin{description}
\item(i) $(\mathcal{A}_{\mathbold{\alpha}}, \mathcal{B}_\mathbold{\alpha})$ is a Cartan pair
\item(ii) there is a homomorphism $\rho_{\mathbold{\alpha}}: \mathcal{A}_\mathbold{\alpha}\to \frac{\prod_n \mathcal{C}_n}{\bigoplus_n \mathcal{C}_n}$ where $\mathcal{C}_n = \mathcal{B}(L^2(X)\otimes l^2\Gamma)$, and $(\rho_{\mathbold{\alpha}}(\mathcal{A}_\mathbold{\alpha}), \rho_{\mathbold{\alpha}}(\mathcal{B}_\mathbold{\alpha}))$ is a Cartan pair
\item(iii) the restrictions of $\rho_{\mathbold{\alpha}}$ to $\mathcal{B}_\mathbold{\alpha}$ and to $\mathcal{N} = C^*\Gamma$, the normalizer algebra of the pair $(\mathcal{A}_\mathbold{\alpha}, \mathcal{B}_{\mathbold{\alpha}})$, have natural lifts, $\psi$ and $\phi$, respectively.
\end{description}}\label{Cartan characterization}

\noindent The dynamical perturbing problem associated to $\mathbold{\alpha}$ has a solution if and only if the Cartan lifting problem associated to $(\mathcal{A}_\mathbold{\alpha}, \mathcal{B}_\mathbold{\alpha}, \rho_{\mathbold{\alpha}}, \psi)$ has a solution.  

\begin{proof}

We can assume WLOG that the $\alpha_n$ respect inverses and that $\alpha_n(e)$ is always the identity.  By \ref{freeness lemma}, we can replace $(\alpha_n)$ by the associated action constructed in that lemma.  By assumption, we have

$$\sup_{\gamma, \lambda\in \Gamma} \sup_{x\in X_n} d_{X_n}(\alpha_n(\gamma)\circ\alpha_n(\lambda)(x), \alpha_n(\gamma\lambda)(x))\to 0.$$

\noindent as $n\to\infty$.  There is a natural (left-regular) action $\Gamma\acts \prod_{\gamma\in \Gamma} C(X_n)$.  Let 

$$B_n = \Gamma\cdot \{(b_\gamma) \mid b_\gamma = \widehat{\alpha}_n(\gamma)(b)\text{, }b\in C(X_n)\}$$  

\noindent and note that this is a $C^*$-algebra.  


Let 
$$\mathcal{B}_{\mathbold{\alpha}} = \{(b_n) \mid b_n\in B_n\text{, and }\{b_n(\gamma_n)\}\text{ is equicontinuous for some sequence }(\gamma_n)_{n\in \mathbb{N}}\}.$$

\noindent Note that $\mathcal{B}_{\mathbold{\alpha}}$ is commutative. If $\mathcal{A}_{\mathbold{\alpha}} = \mathcal{B}_{\mathbold{\alpha}}\rtimes \Gamma$, then $(\mathcal{A}_\mathbold{\alpha}, \mathcal{B}_\mathbold{\alpha})$ is a Cartan pair.  To see why $B_\mathbold{\alpha}$ is maximal abelian in $\mathcal{B}_{\mathbold{\alpha}}\rtimes \Gamma$, let $c_{n, \gamma}\in B_n$ be the function which takes $\lambda$ to the characteristic function associated to $X_{n, \gamma\lambda}$ (the copy of $X_n$ in the disjoint union $\bigsqcup_{\gamma\in \Gamma} X_n$ with index $\gamma\lambda$), and note that $\{c_{n, \gamma} \mid n\in \mathbb{N}, \gamma\in \Gamma\}$ is an equicontinuous family of functions.  Now, if $\sum_{i=1}^{I_n} b_{i, n}u_{\gamma_i}\in B_n\rtimes \Gamma\setminus B_n$ (meaning $b_{i_0, n}\neq 0$ for some $i_0$ with $\gamma_{i_0}\neq e$), there is $\gamma\in \Gamma$ such that 

$$\|(\sum_{i=1}^{I_n} b_{i, n}u_{\gamma_i})c_{n, \gamma} - c_{n, \gamma}(\sum_{i=1}^{I_n} b_{i, n}u_{\gamma_i})\| =\|\sum_{i=1}^{I_n} (b_{i, n} c_{n, \gamma_i\gamma} - b_{i, n}c_{n, \gamma})u_{\gamma_i}\|$$
$$\geq \sup_{i, \lambda}\| b_{i, n}(\lambda)c_{n, \gamma_i\gamma}(\lambda) - b_{i, n}(\lambda)c_{n, \gamma}(\lambda)\|\geq \|b_{i_0, n}\|_\infty>0$$

So $\mathcal{B}_\mathbold{\alpha}$ is maximal abelian in $\mathcal{A}_\mathbold{\alpha}$, implying $(\mathcal{A}_\mathbold{\alpha},\mathcal{B}_\mathbold{\alpha})$ is a Cartan pair. Define $\psi_n((b_\gamma)_{\gamma\in \Gamma})$ to be the operator in $\mathcal{B}(L^2(X_n)\otimes l^2\Gamma)$ which takes $f\otimes \delta_{\gamma_0}\mapsto M_{\widehat{\alpha}_n(\gamma_0^{-1})(b_n(e))}f\otimes \delta_{\gamma_0}$.  Observe that $u_\gamma\psi_n(c_{n, e})u_\gamma^* = \psi_n(c_{n, \gamma^{-1}}) = \psi_n(\gamma\cdot c_{n, e})$. By a slight abuse of notation, we will also use $\lambda$ for the left-regular representation of $\Gamma$. Then if $f\in C(X_n)$, $\delta_{\gamma_0}\in l^2\Gamma$, and $\gamma\in \Gamma$,

$$\psi_n(\gamma\cdot(b_\lambda)_{\lambda\in \Gamma}) = (f\otimes \delta_{\gamma_0}\mapsto M_{\widehat{\alpha}_n(\gamma_0^{-1})(b_{\gamma^{-1}, n})}f\otimes \delta_{\gamma_0})$$

\noindent and

$$(1\otimes \lambda_\gamma)^*\psi_n((b_\lambda)_{\lambda\in \Gamma})(1\otimes \lambda_\gamma) = (f\otimes \delta_{\gamma_0}\mapsto M_{\widehat{\alpha}_n(\gamma_0^{-1}\gamma)(b_n(e))}f\otimes \delta_{\gamma_0})$$

\noindent Furthermore,

$$[\widehat{\alpha}_n(\gamma_0^{-1})(b_n(\gamma^{-1})](x) = [\widehat{\alpha}_n(\gamma_0^{-1})(\widehat{\alpha}_n(\gamma)(b_n(e)))](x)$$
$$= [\widehat{\alpha}_n(\gamma)(b_n(e))](\alpha_n(\gamma_0)(x)) = b_n(e)(\alpha_n(\gamma^{-1})\circ \alpha_n(\gamma_0)(x))$$
$$ = [\widehat{\alpha}_n(\gamma_n)](b_n(e))(\alpha_n(\gamma^{-1})\circ \alpha_n(\gamma_0)\circ \alpha_n(\gamma_n)(x))$$

\noindent and

$$[\widehat{\alpha}_n(\gamma_0^{-1}\gamma)(b_n(e))](x) = b_n(e)(\alpha_n(\gamma^{-1}\gamma_0)(x)) = [\widehat{\alpha}_n(\gamma_n)(b_n(e))](\alpha_n(\gamma^{-1}\gamma_0)\circ \alpha_n(\gamma_n)(x)).$$

\noindent So, if $(b_n(\gamma_n))_{n\in \mathbb{N}} = (\widehat{\alpha}_n(\gamma_n)(b_n(e)))_{n\in \mathbb{N}}$ is equicontinuous, 

$$\|\psi_n(\gamma\cdot(b_n)) - (1\otimes \lambda_\gamma)^*\psi_n(b_n)(1\otimes \lambda_\gamma)\|\to 0$$

\noindent as $n\to\infty$.  That is, the failure of covariance tends to zero as $n\to\infty$.  Since $\mathcal{B}_{\mathbold{\alpha}}$ is generated by elements of the form $\gamma\cdot (b_n)_{n\in \mathbb{N}, \lambda\in \Gamma}$ where $(b_n(\gamma_n)) = ((\widehat{\alpha}_n)(\gamma_n)(b_n(e)))$ is equicontinuous, $(\psi_n)_{n\in \mathbb{N}}$ and $(1\otimes \lambda)$ give a covariant representation, i.e. a homomorphism

$$\rho_{\mathbold{\alpha}}: \mathcal{A}_{\mathbold{\alpha}}\to \frac{\prod_n \mathcal{B}(L^2(X_n)\otimes l^2\Gamma)}{\bigoplus_n \mathcal{B}(L^2(X_n)\otimes l^2\Gamma)}.$$

By our earlier observations about the collection $\{c_{n, \gamma} \mid n\in\mathbb{N}\text{, }\gamma \in \Gamma\}$, if $a\in \mathcal{A}_\mathbold{\alpha}$ and $\rho_\mathbold{\alpha}(a)\notin \rho_{\mathbold{\alpha}}(\mathcal{B}_\mathbold{\alpha})$, then $a = (\sum_{i=1}^{I_n} b_{i, n}u_{\gamma_i})_{n\in \mathbb{N}}$ with $i_n$ and $\epsilon>0$ such that $\|b_{i_n, n}\|\geq \epsilon$ and $u_{\gamma_{i_n}}\neq e$ for all $n$; and so there is a sequence $\gamma_n$ such that $\rho_\mathbold{\alpha}(a)(\psi_n(c_{n, \gamma_n} ))_{n\in \mathbb{N}} - (\psi_n(c_{n, \gamma_n}))_{n\in \mathbb{N}}\rho_{\mathbold{\alpha}}(a)\neq 0$, so $\rho_{\mathbold{\alpha}}(\mathcal{B}_\mathbold{\alpha})$ is maximal abelian in $\rho_{\mathbold{\alpha}}(\mathcal{A}_\mathbold{\alpha})$. Hence, $(\rho_{\mathbold{\alpha}}(\mathcal{A}_\mathbold{\alpha}), \rho_{\mathbold{\alpha}}(\mathcal{B}_\mathbold{\alpha}))$ is a Cartan pair. Let $\mathcal{N}$ be the normalizer algebra associated to the Cartan pair $(\mathcal{A}_{\mathbold{\alpha}}, \mathcal{B}_{\mathbold{\alpha}})$, that is, the algebra $C^*\Gamma$.  The restrictions of $\rho_{\mathbold{\alpha}}$ to $\mathcal{B}_{\mathbold{\alpha}}$ and to $\mathcal{N}$ have natural lifts ($(\psi_n)_{n\in \mathbb{N}}$ and $(1\otimes \lambda)$, respectively).  Also, $\psi_n(\mathcal{B}_{\mathbold{\alpha}})\cong C(X_n)$ for all $n$.

Now we show the perturbing problem associated to $\mathbold{\alpha}$ has a solution if and only if the Cartan lifting problem associated to $(\mathcal{A}_\mathbold{\alpha}, \mathcal{B}_{\mathbold{\alpha}}, \rho_{\mathbold{\alpha}}, \psi)$ has a solution.  

Assume we have a solution to the Cartan lifting problem.  By conditions (1) and (2) of \ref{Cartan lift definition}, $\|\widetilde{\Phi}_n(u_\gamma)^*\widetilde{\Phi}_n(b)\widetilde{\Phi}_n(u_\gamma) - (1\otimes\lambda_\gamma)^*\psi_n(b)(1\otimes \lambda_\gamma)\|\to 0$ as $n\to\infty$ for all $b\in \mathcal{B}_\mathbold{\alpha}$ and $\gamma\in \Gamma$. Observe also that if $\sum_{\gamma\in \Gamma} s_\gamma\lambda_\gamma$ ($s_\gamma\in \mathbb{C}$) is a linear combination, there is $b\in \mathcal{B}_\mathbold{\alpha}$ with $b_n = c_{n, \gamma_n}$ such that 

$$\|\Big(\sum_{\gamma\in \Gamma} \bar{s}_\gamma(1\otimes \lambda^*_\gamma)\Big)\psi(b)\Big(\sum_{\gamma\in \Gamma} s_\gamma(1\otimes \lambda_\gamma)\Big)\|\geq \max_{\gamma\neq \gamma'} |\bar{s}_\gamma s_{\gamma'}|.$$ 

\noindent So any non-trivial linear combination is far from normalizing $\psi(b)$ for some such $b$. Hence, the same is true for any non-trivial linear combination in $\{\widetilde{\Phi}(u_\gamma) \mid \gamma\in \Gamma\}$. But by condition (3) of \ref{Cartan lift definition}, the normalizers of $\psi(\mathcal{B}_\mathbold{\alpha})$ must generate $\widetilde{\Phi}(\mathcal{A}_\mathbold{\alpha})$, so $\widetilde{\Phi}(u_\gamma)$ normalizes $\widetilde{\Phi}(\mathcal{B}_\mathbold{\alpha}) = \psi(\mathcal{B}_\mathbold{\alpha}) \cong C(X_n)$ for all $\gamma\in \Gamma$.  So each $\widetilde{\Phi}(u_\gamma)$ corresponds to a homeomorphism $X_n\to X_n$.  By condition (2), $u_\gamma\to \widetilde{\Phi}(u_\gamma)$ is an approximate representation of $\Gamma$, and so these homeomorphisms $X_n\to X_n$ are uniformly close to giving a group action.  But the only way that can happen is if they actually do give a group action, and so we have actions on $\psi_n(\mathcal{B}_{\mathbold{\alpha}})\cong C(X_n)$ and therefore actions on $X_n$ by duality.  Since $\pi\circ \Phi: \mathcal{A}_\mathbold{\alpha}\to \frac{\prod_n \mathcal{C}_n}{\bigoplus_n \mathcal{C}_n}$ is a homomorphism, the maps $\mathcal{B}_{\mathbold{\alpha}}\xrightarrow{\psi_n} C(X_n)$ are approximately equivariant in the sense that, for a fixed $b\in \mathcal{B}_{\mathbold{\alpha}}$, the failure of equivariance between $\Gamma\acts \mathcal{B}_\mathbold{\alpha}$ and these new actions on $C(X_n)$ is eventually less than $\epsilon$ for large $n$.  Since each $C(X_n)$ is generated by at most $M$ functions, and the union of all these generating sets is equicontinuous, we have a finite subset of $\mathcal{B}_{\mathbold{\alpha}}$ whose image in each $C(X_n)$ contains a generating set.  Hence, for all $\epsilon$ there is $N$ such that, for $n\geq N$, the aforementioned map is within $\epsilon$ of being equivariant on a generating set of $C(X_n)$.  The newly-constructed actions are therefore perturbations of the original almost-actions.  

Now suppose we have a solution, $(\tilde{\alpha}_n)$, to the perturbing problem associated to $(\alpha_n)$.  By slight abuse of notation, we will use also $\widetilde{\alpha}_n$ for the dual actions on $C(X_n)$. Construct $\mathcal{B}_\mathbold{\widetilde{\alpha}}$ and $\widetilde{\psi}_n$ the same way as $\mathcal{B}_\mathbold{\alpha}$ and $\psi_n$ but with $(\widetilde{\alpha}_n)$ in place of $(\alpha_n)$.  Notice that $(\widetilde{\psi}_n)$ lifts $(\psi_n)$ and that $\widetilde{\psi}_n(\mathcal{B}_\mathbold{\alpha})\cong C(X_n)\cong \psi_n(\mathcal{B}_\mathbold{\alpha})$. Then $(\widetilde{\psi}_n)$ lifts $(\psi_n)$ is approximately covariant with $1\otimes \lambda$, so if we use these two maps to define a $*$-linear map $\widetilde{\Phi}: \mathcal{A}_\mathbold{\alpha}\to \prod_n \mathcal{C}_n$, then $\widetilde{\Phi}$ gives a homomorphism $\mathcal{A}_\mathbold{\alpha}\to \frac{\prod_n\mathcal{C}_n}{\bigoplus_n \mathcal{C}_n}$.  Now

$$u_\gamma\widetilde{\psi}_n(b_n)u_\gamma^* = (f\otimes \delta_\lambda\mapsto M_{\widetilde{\alpha}_n(\lambda\gamma^{-1})(b_{n}(e))} f\otimes \delta_\lambda)$$
$$=(f\otimes \delta_\lambda\mapsto M_{\widetilde{\alpha}_n(\lambda)\circ \widetilde{\alpha}_n(\gamma^{-1})(b_{n}(e))}f\otimes \delta_\lambda)$$

\noindent and $\widetilde{\alpha}_n(\lambda)\circ \widetilde{\alpha}_n(\gamma^{-1})(b_{n}(e))$ is equicontinuous for some sequence $\lambda_n\in \Gamma$. Hence, $\widetilde{\Phi}(u_\gamma)$ normalizes $\widetilde{\Phi}(b)$ for $b\in \mathcal{B}_\mathbold{\alpha}$. 

Put $(d_n(e)) = (\widetilde{\alpha}_n(\gamma^{-1})(b_n(e))$.  Then if $(b_n(\gamma_n))$ is equicontinuous,

$$\widehat{\alpha}_n(\lambda)\circ \widetilde{\alpha}_n(\gamma^{-1})(b_n(e)) = \widehat{\alpha}_n(\lambda)\circ \widetilde{\alpha}_n(\gamma^{-1})\circ \widehat{\alpha}_n(\gamma_n^{-1})(b_n(\gamma_n))$$
$$\approx_\epsilon \widehat{\alpha}_n(\lambda)\circ \widetilde{\alpha}_n(\gamma^{-1})\circ \widetilde{\alpha}_n(\gamma_n^{-1})(b_n(\gamma_n)) = \widehat{\alpha}_n(\lambda)\circ \widetilde{\alpha}_n(\gamma^{-1}\gamma_n^{-1})(b_n(\gamma_n))$$
$$\approx_\epsilon \widehat{\alpha}_n(\lambda)\circ \widehat{\alpha}_n(\gamma^{-1}\gamma_n^{-1})(b_n(\gamma_n)) \approx_\epsilon \widehat{\alpha}_n(\lambda\gamma^{-1}\gamma_n^{-1})(b_n(\gamma_n)) = b_n(\gamma_n)$$

\noindent when $\lambda = \gamma_n\gamma$. So there are $\lambda_n\in \Gamma$ such that $(\widehat{\alpha}_n(\lambda)\circ \widetilde{\alpha}_n(\gamma^{-1})(b_n(e)))$ is equicontinuous. Notice that being of the form $\gamma\mapsto \widehat{\alpha}_n(\lambda)(b_n(e))$ with $\alpha_n(\lambda_n)(b_n(e))$ equicontinuous for some $\lambda_n$ is equivalent to being of the form $\lambda_n\cdot b_n$ with $(b_n(e))$ equicontinuous. So if $(b_n) = b\in \mathcal{B}_\mathbold{\alpha}$, there is $d = (d_n)$ such that $\widetilde{\psi}(d) = u_\gamma\widetilde{\psi}(b)u_\gamma^* = \gamma\cdot\widetilde{\psi}(b)$, where at the end, we're thinking of $\widetilde{\psi}(b)$ as an element of $\mathcal{B}_{\widetilde{\mathbold{\alpha}}}$. Hence, $\widetilde{\psi}(\mathcal{B}_\mathbold{\alpha}) = \mathcal{B}_{\widetilde{\mathbold{\alpha}}}$. Since $(\widetilde{\Phi}(\mathcal{A}_{{\mathbold{\alpha}}}), \widetilde{\Phi}(\mathcal{B}_{{\mathbold{\alpha}}}))$ is therefore isomorphic to $(\mathcal{A}_{\widetilde{\mathbold{\alpha}}}, \mathcal{B}_{\widetilde{\mathbold{\alpha}}})$, and since the latter is a Cartan pair by the same argument we used to show $(\mathcal{A}_{\mathbold{\alpha}}, \mathcal{B}_{\mathbold{\alpha}})$ is (using that the actions $\widetilde{\alpha}_n$ are free), $(\widetilde{\Phi}(\mathcal{A}_\mathbold{\alpha}), \widetilde{\Phi}(\mathcal{B}_\mathbold{\alpha}))$ is a Cartan pair.


Observe that there is an inner automorphism of $\prod_n \mathcal{C}_n$ which restricts to an isomorphism $\text{im}(\widetilde{\psi}_n) \cong \text{im}(\psi_n)$. Replace $1\otimes \lambda$ by its conjugation by this automorphism) and replace $\widetilde{\Phi}$ with the $*$-linear map given rise to by this new representation and $\psi = (\psi_n)$. This new $\widetilde{\Phi}$ is then a solution to the Cartan lifting problem. \end{proof}

%
%
%
%

\remark{It is worth remarking that one could prove a similar theorem if, in \ref{Cartan lift definition}, we replaced (1) with the requirements that $\widetilde{\Phi}_{|\mathcal{B}_\mathbold{\alpha}}$ lifts the restriction of $\rho_\mathbold{\alpha}$ to $\mathcal{B}$, and has image isomorphic to $\psi(\mathcal{B})$, and that $\widetilde{\Phi}_{|C^*\Gamma} = 1\otimes \lambda$.  This resolves some of the apparent disanology between these problems and traditional lifting problems. While it might appear more natural in the definition to require that $\widetilde{\Phi}_{|C^*\Gamma}$ and $\widetilde{\Phi}_{|\mathcal{B}_\mathbold{\alpha}}$ be lifts, the proof above shows the condition given in \ref{Cartan lift definition} is actually equivalent to something apparently stronger, so we are in fact asking about the existence of some sort of special lifts. Rather than having the ``closeness" condition enforced by explicitly requiring a map to be a lift, one can view conditions (1) and (2) of \ref{Cartan lift definition} as enforcing the closeness condition. }\label{alt definition remark}

\noindent \normalfont An elementary argument also shows that almost-actions by finite groups on Cantor sets are always near actual actions. 

\theorem{Suppose $\Gamma$ is a finite group, $X$ is a Cantor set, and $\alpha_n: \Gamma\to \text{Homeo}(X)$ is an almost-action. The dynamical perturbing problem associated to $(\alpha_n)$ has a solution.} \label{pretty Cantor sets theorem}
\begin{proof}
Given $\epsilon>0$, find a finite, clopen partition $\{X_i\}_{i=1}^I$ of $X$ into sets with diameter $<\epsilon$.  Let $\delta>0$ be the least distance between two subsets in this partition.  Let $\mathcal{X}$ be the join of the collections $\{\gamma\cdot X_i\}_{i=1}^I$ for all $\gamma\in \Gamma$, that is, the collection of all sets of the form $\bigcap_{\gamma\in \Gamma} \gamma\cdot X_{i_\gamma}$.  Find $N$ such that $d_X(\alpha_n(\gamma)\circ \alpha_n(\delta)(x), \alpha_n(\gamma\delta)(x))<\delta$ for all $\gamma, \delta\in \Gamma$ for $n\geq N$.  Then for such $n$, $\alpha_n$ preserves $\mathcal{X}$.  Since any two clopen subsets of $X$ are homeomorphic, we can define an (actual) action $\Gamma\acts X$ which agrees with the action on $\mathcal{X}$ given by $\alpha_n$, and therefore stays within $\epsilon$ of $\alpha_n$.  We therefore have a solution to the perturbing problem associated to $(\alpha_n)$. \end{proof}

\normalfont Finally, we can find solutions when $\Gamma$ is finite and the almost-actions are equicontinuous. We include this last result mainly to further establish that there are many non-trivial examples where perturbing problems (and therefore associated Cartan lifting problems) have solutions. 

\theorem{Suppose $\Gamma$ is a finite group, $X$ is a compact metric space, and $\alpha_n: \Gamma\to \text{Homeo}(X)$ are set maps such that $\{\alpha_n(\Gamma)\}$ is equicontinuous independent of $n$, and 

$$d_X(\alpha_n(\gamma_1)\circ \alpha_n(\gamma_2)(x), \alpha_n(\gamma_1\gamma_2)(x))\to 0\text{ as }n\to\infty$$
for all $\gamma_1, \gamma_2\in \Gamma$ and $x\in X$.  Then there are a subsequence $(n_k)$ and an action $\tilde{\alpha}: \Gamma\to \text{Homeo}(X)$ such that $\sup_{x\in X}d_X(\tilde{\alpha}(\gamma)(x), \alpha_{n_k}(\gamma)(x))\to 0$ as $k\to\infty$ for all $\gamma\in \Gamma$.  }

\begin{proof}
Let $S=\{\alpha_n(\gamma) \mid \gamma\in \Gamma\text{, }n\in \mathbb{N}\}$.  Then $S$ is an equicontinuous family of homeomorphisms $X\to X$, and so the closure of $S$ is compact by the Arzel\`{a}-Ascoli theorem.  By repeatedly extracting nested subsequences, we can therefore find a subsequence $n_k$ such that, for all $\gamma\in \Gamma$, $\alpha_{n_k}(\gamma)$ converges to some $\beta(\gamma)\in \text{Homeo}(X)$.  It can be checked that $\gamma\to \beta(\gamma)$ is a group action.  \end{proof}

\normalfont
\noindent In light of these results, we conclude with the following conjecture. 

\conjecture{If $\Gamma$ is a finite group, $X$ is a compact metric space, and $\alpha_n:\Gamma\to X$ is an almost-action, the perturbing problem associated to $(\alpha_n)$ has a solution. }

\normalfont While it seems plausible that the semiprojectivity of finite-dimensional $C^*$-algebras, or the $C^*$-algebraic characterization of perturbing problems given above may provide a means to prove this conjecture, it is not obvious how, at least to the author at present.

\bibliography{mybibliography2.bib}
\bibliographystyle{plain}

\end{document}